\documentclass[11pt]{amsart}
\renewcommand{\and}{et}

\usepackage{array}
\usepackage{amssymb,amsthm,amsfonts,amssymb,euscript,mathrsfs
}
\usepackage[english,francais]{babel}
\usepackage{amsxtra}
\usepackage{amsrefs}
\usepackage{amscd}
\usepackage{enumitem}
\usepackage[T1]{fontenc}                 
\usepackage{color}

 \DeclareMathAlphabet{\got}{U}{euf}{m}{n}     
\DeclareMathAlphabet{\mat}{U}{msb}{m}{n}     
\DeclareMathAlphabet{\mathbold}{OML}{cmm}{bx}{it} 

\makeatletter

\renewcommand{\section}{\@startsection{section}{1}%
\z@{-3.5ex \@plus -1ex \@minus -.2ex}
  {2.3ex \@plus.2ex}
  {\bfseries\large}}           

\renewcommand{\subsection}{\@startsection{subsection}{2}%
\z@{.5ex \@plus.7ex}{-.5em}%
{\normalfont\bfseries}}

\renewcommand{\subsubsection}{\@startsection{subsubsection}{3}%
\z@{.5ex\@plus.7ex}{-.5em}%
  {\normalfont\itshape}}
\makeatother

\makeatletter
\renewcommand{\paragraph}{\@startsection{paragraph}{4}
  \z@{.5ex\@plus.7ex}{-.5em}%
                                    {\normalfont\normalsize\itshape}}%
\makeatother
\let\subsubsubsection\paragraph
\numberwithin{equation}{section}

\newtheorem{theo}{Th\'eor\`eme}[]
\newtheorem{pr}{Proposition}[]
\newtheorem{lem}{Lemme}[]

\newenvironment{rem}{\vskip 1em\bf Remarque.\rm}{\par\rm}
\newenvironment{theo-1}{\vskip 1em\bf Theorem\it}{\par\rm}
\newenvironment{theo*}{\vskip 1em\bf Th\'eor\`eme.\it}{\par\rm}
\newenvironment{dem}{\vskip 1em{\it D\'emonstration} :}%
{\unskip\hfill\null\nobreak\hfill\carre\vskip1em\par}
\newcommand{\carre}{\rule{1ex}{1ex}}

\makeatletter
\providecommand*{\diff}%
{\@ifnextchar^{\DIfF}{\DIfF^{}}}
\def\DIfF^#1{%
\mathop{\mathrm{\mathstrut d}}%
\nolimits^{#1}\gobblespace}
\def\gobblespace{%
\futurelet\diffarg\opspace}
\def\opspace{%
\let\DiffSpace\!%
\ifx\diffarg(%
\let\DiffSpace\relax
\else
\ifx\diffarg[%
\let\DiffSpace\relax
\else
\ifx\diffarg\{%
\let\DiffSpace\relax
\fi\fi\fi\DiffSpace}

\providecommand*{\eu}%
{\ensuremath{\mathrm{e}}}

\providecommand*{\iu}%
{\ensuremath{\mathrm{i}}}

\newcommand{\Ad}{\mathop{\mathrm{Ad}}}
\newcommand{\im}{\mathop{\mathrm{Im}}}

\newcommand{\ad}{\mathop{\mathrm{ad}}}

\newcommand{\Hom}{\mathop{\mathrm{Hom}}}

\newcommand{\rg}{\mathop{\mathrm{rg}}}
\newcommand{\diag}{\mathop{\mathrm{diag}}}

\newcommand{\Tr}{\mathop{\mathrm{Tr}}}

\begin{document}

\title{Discriminants et Sommes de carrés}

\date{\today}
\selectlanguage{french}
\author{Mustapha Raïs}
\address{Poitiers}
\email{mustapaha.rais@math.univ-poitiers.fr}
\thanks{L'auteur remercie chaleureusement
    Pierre Torasso pour sa généreuse contribution à l'élaboration de
    ce texte}





\maketitle

\hfill\parbox[b]{5cm}{En hommage à Claude Quitté pour son grand
  dévouement à l'enseignement et à la recherche mathématiques}

\section{Introduction}\label{1}
\subsection{}\label{1.1}
On note $\mathfrak{p}$ l'espace vectoriel des matrices symétriques
réelles de taille $n$ et $\mathfrak{k}=so(n)$ celui des matrices
réelles antisymétriques de même taille. Lorsque $x\in\mathfrak{p}$, le
discriminant de la matrice $x$ se calcule de plusieurs manières, par
exemple :
\begin{equation*}
D(x)=\prod_{i<j}(\lambda_{i}-\lambda_{j})^{2}\mbox{, les $\lambda_{k}$
  étant les valeurs propres de $x$}
\end{equation*}
Par ailleurs, la fonction $D:\mathfrak{p}\rightarrow\mathbb{R}$ est
une fonction polynôme qui s'exprime par l'intermédiaire des
coefficients du polynôme caractéristique de $x$, puisque c'est le
résultant du polynôme caractéristique et de sa dérivée ; c'est donc la
seule fonction $K$-invariante sur $\mathfrak{p}$, dont la valeur sur
les matrices $\diag(\lambda_{1},\ldots,\lambda_{n})$ est
$\prod_{i<j}(\lambda_{i}-\lambda_{j})^{2}$. Ici $K=SO(n)$ est le
groupe des rotations et l'opération de $K$ dans $\mathfrak{p}$ est
l'opération de conjugaison.

\subsection{}\label{1.2}
Il y a d'autres expressions de $D$ comme fonction polynôme sur
$\mathfrak{p}$ (en dehors de celle classique donnant le discriminant
d'un polynôme sous forme d'un déterminant où apparaissent les
coefficients du polynôme (et de sa dérivée)). Pour présenter la
première, on note $\langle\;\vert\;\rangle$ le produit scalaire 
\begin{equation*}
\langle x\vert y\rangle=\Tr(x\;^{t}\!y)
\end{equation*}
sur $M_{n}(\mathbb{R})$ (celui pour lequel la base canonique
$(E_{ij})_{i,j}$ est une base orthonormée), et
$p_{1}(x),$ $p_{2}(x),\ldots,p_{n}(x)$ les coefficients du polynôme
caractéristique de $x$ :
\begin{equation*}
  \det(tI_{n}-x)=t^{n}-p_{1}(x)t^{n-1}+\cdots+(-1)^{n}p_{n}(x)
\end{equation*}
Chaque $p_{j}:\mathfrak{p}\rightarrow\mathbb{R}$ admet un grandient
$\nabla p_{j}:\mathfrak{p}\rightarrow\mathfrak{p}$ calculé au moyen de
la restriction à $\mathfrak{p}$ du produit scalaire ci-dessus. On
introduit alors la matrice :
\begin{equation*}
\Omega(x)=(\Omega_{ij}(x))\mbox{, où $\Omega_{ij}(x)=\langle\nabla
  p_{i}(x)\vert\nabla p_{j}(x)\rangle$}
\end{equation*}
et il vient :
\begin{equation}\label{eq1.1}
D(x)=\det \Omega(x)
\end{equation}

\subsection{}\label{1.3}
Pour deux autres présentations de $D$, on introduit pour chaque $x$
dans $\mathfrak{p}$, quatre applications linéaires :
\begin{equation*}
A(x) : \mathfrak{p}\rightarrow\mathfrak{k}\mbox{, }B(x) :
\mathfrak{k}\rightarrow\mathfrak{p}
\end{equation*}
(chacune ici transforme $y$ en $[x,y]=(\ad x).y$)
\begin{align*}
  f(x)&=B(x)\circ A(x) : \mathfrak{p}\rightarrow\mathfrak{p}\\
  g(x)&=A(x)\circ B(x) : \mathfrak{k}\rightarrow\mathfrak{k}
\end{align*}
(chacune ici envoie $y$ sur $[x,[x,y]]=(\ad x)^{2}.y$).
On écrit les polynômes caractéristiques respectifs de $f(x)$ et $g(x)$
:
\begin{align*}
\det(tId_{\mathfrak{p}}-f(x))&=t^{r}-P_{1}(x)t^{r-1}+\cdots+(-1)^{r-l}P_{r-l}(x)t^{l}\\
\det(tId_{\mathfrak{k}}-g(x))&=t^{s}-Q_{1}(x)t^{s-1}+\cdots+(-1)^{s}Q_{s}(x)  
\end{align*}
avec $r=\frac{n(n+1)}{2}$, $l=n$, $s=\frac{n(n-1)}{2}$ et il vient :
\begin{equation}\label{eq1.2}
  D(x)=P_{r-l}(x)=Q_{s}(x)
\end{equation}
A ce point, des remarques :
\begin{itemize}[label=-]
\item $g(x)$ est \flqq\;génériquement bijective \frqq~ et
\begin{equation*}
D(x)=\det g(x)
\end{equation*}

\item On a : 
\begin{equation*}
  t^{r}\det(tId_{\mathfrak{k}}-g(x))=t^{s}\det(tId_{\mathfrak{p}}-f(x))
\end{equation*}

\item Le produit scalaire $\langle\;\vert\;\rangle$ fait de
$\mathfrak{k}$ et $\mathfrak{p}$ des espaces euclidiens et $B(x)$
(resp. $A(x)$) est l'adjoint de $A(x)$ (resp. $B(x)$) :
$B(x)=A(x)^{*}$, $A(x)=B(x)^{*}$.
\end{itemize}

\subsection{}\label{1.4}
Certaines fonctions polynômes $f(x_{1},\ldots,x_{n})$
($f\in\mathbb{R}[x_{1},\ldots,x_{n}]$) sont des sommes de carrés :
\begin{equation*}
  f=\sum_{1\leq j\leq k}g_{j}^{2}\mbox{\; $g_{j}\in\mathbb{R}[x_{1},\ldots,x_{n}]$}
\end{equation*}
Il en est ainsi lorsque, par exemple, $f$ est une forme quadratique
vérifiant la condition nécessaire de positivité : $f(x)\geq0$ pour
tout $x$, mais il existe des fonctions polynômes positives qui ne sont
pas des sommes de carrés (voir \cite{prestel-delzell-b-2001} et
\cite{rudin-2000}).

Par contre, le discriminant des matrices symétriques est une somme de
carrés.

D'après M\'{a}ty\'{a}s Domokos, ce résultat est dû à Borchardt (1846)
; voir \cite{domokos-2011} qui est pour une grande part à l'origine de
ce texte.

Ceci étant \flqq\;être une somme de carrés \frqq~ peut donner lieu à
diverses variations. D'abord, lorsque $f$ est une somme de carrés, on
peut chercher à déterminer le nombre minimum de polynômes $g_{j}$ qui
interviennent dans les diverses décompositions :
\begin{equation*}
f=\sum g_{j}^{2}
\end{equation*}
C'est pour l'essentiel cette question pour le discriminant $D$ qui
motive le travail de Domokos. On trouvera plus loin des indications
sur le travail de Domokos et aussi sur ceux d'autres
mathématiciens. Ensuite, on peut restreindre le choix des $g_{j}$ à
une classe particulière de polynômes. Dans cet ordre d'idée, il est
indiqué d'exhumer un intéressant article de Walter Rudin \cite{rudin-1984}.

\subsection{}\label{1.5}
Les fonctions $f$ auxquelles Rudin s'intéresse sont les puissances des
normes hermitiennes (resp. euclidiennes) sur $\mathbb{C}^{n}$
(resp. $\mathbb{R}^{n}$) : $(\sum\vert z_{j}\vert^{2})^{k}$
(resp. $(\sum x_{j}^{2})^{k}$) qu'il cherche à écrire comme somme de
carrés à l'exemple de : $(\vert z_{1}\vert^{2}+\vert
z_{2}\vert^{2})^{2}=\vert g_{1}(z)\vert^{2}+\vert
  g_{2}(z)\vert^{2}+\vert g_{3}(z)\vert^{2}$ avec
  $g_{1}(z)=z_{1}^{2}$, $g_{2}(z)=\sqrt{2}z_{1}z_{2}$,
  $g_{3}(z)=z_{2}^{2}$, qui fournit une application holomorphe (et
  même polynômiale) propre :
  $(z_{1},z_{2})\mapsto(g_{1}(z),g_{2}(z),g_{3}(z))$ de
  $\mathbb{C}^{2}$ dans $\mathbb{C}^{3}$ (exemple d'Alexander). Pour
  celà, Rudin étudie la classe $W(k,n,p)$ des applications
  $\phi=(g_{1},\ldots,g_{p})$ de $\mathbb{C}^{n}$ dans
  $\mathbb{C}^{p}$ telles que :
\begin{equation*}
(\sum_{1\leq j\leq n}\vert z_{j}\vert^{2})^{k}=\sum_{1\leq j\leq
    p}\vert g_{j}(z)\vert^{2}
\end{equation*}
les $g_{j}$ étant des polynômes homogènes (de degré $k$) en
$z_{1},\ldots,z_{n}$, et linéairement indépendants. De façon
miraculeuse, Rudin obtient le résultat suivant (où $N$ est la
dimension de l'espace vectoriel $P(k,n)$ des polynômes homogènes de
degré $k$ en $z_{1},\ldots,z_{n}$) :
\begin{itemize}[label=-]
\item Il n'y a pas de solution (i.e. $W(k,n,p)$ est vide) lorsque
  $p\neq N$

\item Lorsque $p=N$, les solutions sont celles $(g_{1},\ldots,g_{p})$
  où $(\sqrt{N}g_{1},\ldots,\sqrt{N}g_{p})$ est une base orthonormée
  de $P(k,n)$.
\end{itemize}
  
Ce qui n'apparaît pas dans cet énoncé, mais qui apparaît clairement
dans la démonstration qu'en donne Rudin, c'est qu'il y a une opération
\flqq\;naturelle \frqq~ du groupe unitaire $U(n)$ dans l'espace
$P(k,n)$, qui définit une représentation unitaire irréductible du
groupe unitaire. On voit apparaître ainsi une intervention de
l'analyse harmonique d'un groupe de Lie dans une question de
décomposition en somme de carrés, et il semble bien que cette
antériorité n'ait pas été notée.

Plus précisément, il existe des réels $c_{\alpha}>0$ tels que
$(c_{\alpha}z^{\alpha})_{\alpha}$,
$\alpha=(\alpha_{1},\ldots,\alpha_{n})$,
$z^{\alpha}=z_{1}^{\alpha_{1}}\cdots z_{n}^{\alpha_{n}}$, soit une bas
orthonormée de $P(k,n)$. Il n'et pas trop tard pour préciser que le
produit scalaire sur $P(k,n)$ s'obtient en identifiant $P(k,n)$ à un
sous-espace vectoriel de $L^{2}(\Sigma)$, $\Sigma$ étant la sphère unité de
$\mathbb{C}^{n}$, munie de la mesure invariante de masse totale
$1$. On a donc :
\begin{equation*}
\frac{1}{N}\sum_{\vert\alpha\vert=k}c_{\alpha}^{2}\vert
z^{\alpha}\vert^{2}
=(\vert z_{1}\vert^{2}+\cdots+\vert z_{n}\vert^{2})^{k}
\end{equation*}
Du coup, il vient avec $(x_{1},\ldots,x_{n})$ dans $\mathbb{R}^{n}$ :
\begin{equation*}
\frac{1}{N}\sum_{\vert\alpha\vert=k}c_{\alpha}^{2}
(x^{\alpha})^{2}
=(x_{1}^{2}+\cdots+x_{n}^{2})^{k}
\end{equation*}
ce qui donne une solution au problème analogue sur $\mathbb{R}$. Rudin
indique que le cas réel diffère du cas complexe en ce qu'il y a
d'autres solutions que celles données par les bases orthonormées de
$P(k,n)$ ; il en donne d'ailleurs des exemples. A ce sujet, voici trois
remarques qui n'apparaissent pas dans l'article de Rudin :
\makeatletter
\renewcommand{\theenumi}{\alph{enumi}}
\renewcommand{\labelenumi}{\theenumi.}
\makeatother
\begin{enumerate}
\item Supposons qu'on ait :
\begin{equation*}
\sum_{1\leq j\leq p}f_{j}^{2}(x)=(x_{1}^{2}+\cdots+x_{n}^{2})^{k}
\end{equation*}
Alors : $\sum_{j}(x_{1}^{2}+\cdots+x_{n}^{2})^{2}
f_{j}^{2}(x)=
(x_{1}^{2}+\cdots+x_{n}^{2})^{k+2}$, ce qui montre : $W(k+2,n,p)$ est
non vide dès que $W(k,n,p)$ l'est aussi.\label{rem1.a}

\item Une autre différence est que le groupe $SO(n)$ opère dans
  $P(k,n)$, mais non de manière irréductible. Comme on va le voir, la
  donnée d'un sous-$SO(n)$-module $V$ de $P(k,n)$ fournit une
  décomposition de $(x_{1}^{2}+\cdots+x_{n}^{2})^{k}$ en une somme de
  $q$ carrés, $q$ étant la dimension de l'espace vectoriel $V$. En
  effet, si $\varphi_{1},\ldots,\varphi_{q}$ est une base orthonormée
  de $V$, le noyau reproduisant $\mathcal{K}$ de $V$ est :
  $\mathcal{K}(x,y)=\sum_{1\leq j\leq q}\varphi_{j}(x)\varphi_{j}(y)$,
  et c'est une fonction $SO(n)$-invariante au sens suivant :
\begin{equation*}
\mathcal{K}(g.x,g.y)=\mathcal{K}(x,y)\mbox{, }g\in SO(n)
\end{equation*}
D'où $\mathcal{K}(x,x)=\sum\varphi_{j}(x)^{2}=c(x_{1}^{2}+\cdots+x_{n}^{2})^{k}
$ ($c>0$) car $\mathcal{K}(x,x)$ est une fonction polynôme
$SO(n)$-invariante, homogène de degré $2k$. Sous des formes à peu près
identiques, ce résultat se retrouve dans l'article de Rudin, dans
celui de Peter Lax dont il sera question plus loin, et apparaît dans
des travaux ultérieurs.

Un exemple de sous-espace $V$ est celui $H(k,n)$ des polynômes
harmoniques de degré $k$. On pourra comparer les dimensions
respectives de $H(k,n)$ et $P(k,n)$ :
\begin{align*}
  \dim H(k,n)&=\frac{(2k+n-2)(n+k-3)!}{(n-2)!k!}\\
  \dim P(k,n)&=\frac{(k+n-1)!}{k!(n-1)!}
\end{align*}

\item
\begin{itemize}[label=-]
\item En utilisant la remarque \ref{rem1.a} ci-dessus par exemple, on
  voit que $(x_{1}^{2}+x_{2}^{2})^{2m+1}$ est somme de $2$ carrés :
\begin{equation*}
  (x_{1}^{2}+x_{2}^{2})^{2m+1}=(x_{1}(x_{1}^{2}+x_{2}^{2})^{m})^{2}+
  (x_{2}(x_{1}^{2}+x_{2}^{2})^{m})^{2}
\end{equation*}
(tandis que $(x_{1}^{2}+x_{2}^{2})^{2m}$ est lui-même un carré).

\item Enfin : Rudin signale que la solution complète de la
  détermination des $W(k,n,p)$ dans le cas réel reste à trouver.
\end{itemize}
\end{enumerate}

\subsection{}\label{1.6}
Soit $Z$ une matrice carrée $n\times n$ à coefficients complexes. Dans
\cite{ilyushechkin-1992}, Ilyushechkin définit une matrice $Z_{*}$ de
type $n^{2}\times n$, de la manière suivante : les lignes de $Z_{*}$
sont indexées par les couples $(i,j)$, $1\leq i,j\leq n$, et les
colonnes par les entiers $k$, $0\leq k\leq n$ ; l'élément de $Z_{*}$
se trouvant dans la colonne $k$ et la ligne $(i,j)$ est :
\begin{equation*}
(Z_{*})_{(i,j),k}=\Tr(E_{ji}Z^{k})
\end{equation*}
i.e. le coefficient d'indice $(i,j)$ de la matrice $Z^{k}$.

Le résultat principal est le suivant :
\selectlanguage{english}
\begin{theo-1}
Let $D(Z)$ be the discriminant of the characteristic polynomial of the
normal matrix $Z$ over the complex numbers. Then :
\begin{equation}\label{eq1.3}
\vert D(Z)\vert=\sum_{\alpha}\vert M_{\alpha}(Z)\vert^{2}
\end{equation}
where the summation on the righthand side is taken over all minors of
maximal order of $Z_{*}$.
\end{theo-1}
\selectlanguage{french}
On notera que dans cette égalité apparaissent les modules des nombres
complexes concernés.
  
La démonstration est très jolie : Ilyushechkin montre que les deux
membres de l'égalité \eqref{eq1.3} sont invariants sous l'action
naturelle du groupe $U(n)$ dans l'espace des matrices $n\times n$ (par
conjugaison). Il est clair que c'est surtout l'invariance du deuxième
membre qui n'est pas évidente. {\color{black} Un grand pas est fait si on
  réalise que les $M_{\alpha}(Z)$ sont les coefficients du $n$-vecteur
  :
\begin{equation*}
I_{n}\wedge Z\wedge Z^{2}\wedge\cdots\wedge Z^{n-1}
\end{equation*}
élément de $\wedge^{n}gl(n,\mathbb{C})$ dans la base
$E_{\alpha}=E_{i_{0},j_{0}}\wedge\cdots\wedge E_{i_{n-1},j_{n-1}}$
avec $\alpha$ parcourant les $n$-uplets
$((i_{0},j_{0}),\ldots,(i_{n-1},j_{n-1}))$ de couples $(i,j)$, $1\leq
i,j\leq n$, deux à deux distincts et rangés dans l'ordre
lexicographique (par exemple). Dès lors, le second membre est la norme
de ce $n$-vecteur pour le produit scalaire hermitien faisant de la
base $(E_{\alpha})_{\alpha}$ une base orthonormée, produit scalaire
$U(n)$-invariant.} Il suffit alors de démontrer l'égalité lorsque $Z$
est une matrice diagonale. Les détails manquants sont lisibles dans
l'article cité. On verra plus loin, dans la partie consacrée à ma
modeste contribution au sujet, que les idées et méthodes
d'Ilyushechkin y jouent un rôle essentiel. Il faut donc lui rendre
l'hommage qu'il mérite. Ceci étant, il y a dans son travail un
complément facile concernant les matrices symétriques réelles.

Le discriminant $D(X)$ d'une matrice réelle symétrique $X$ est la
somme des carrés des mineurs d'ordre $n$ de la matrice $X_{*}$. Le
nombre de ces mineurs est la dimension de l'espace vectoriel
$\wedge^{n}\mathfrak{p}$ (où $\mathfrak{p}$ est celui des matrices
symétriques), à savoir le coefficient binomial $
\begin{pmatrix}
  \frac{n(n+1)}{2}\\
  n
\end{pmatrix}$.

\subsection{}\label{1.7}
L'article \cite{lax-1998} de Peter D. Lax est à l'origine de divers
travaux les plus récents portant sur la question suivante : Quel est
le nombre minimum de carrés qui interviennent dans une somme de carrés
égale à $D$? Des exemples portant sur les petites valeurs de $n$ sont
donnés dans l'article de Domokos cité. La contribution de Lax est
originale et remarquable.

Soit $X$ l'ensemble des $x$ dans $\mathfrak{p}$ tels que $D(x)=0$,
i.e. celui des matrices ayant (au moins) une valeur propre
multiple. Soit $V$ l'espace vectoriel des fonctions polynômes dans
$\mathbb{R}(\mathfrak{p})$ qui s'annulent sur $X$ et qui sont
homogènes de degré $\frac{n(n-1)}{2}$ (il n'est pas clair a priori que $V$ 
contienne des polynômes non nuls, et Lax exhibe un élément non nul de
$V$). Soit $(r_{i})_{i}$ une base de $V$. Alors :
\begin{itemize}[label= ]
\item \flqq\;$D$ est une somme de carrés de combinaisons linéaires des
  $r_{i}$ \frqq 
\end{itemize}
Dans cet article, on voit apparaître l'idéal $\mathcal{J}(X)$ de l'ensemble
algébrique $X$ et l'action de $SO(n)$ sur $\mathfrak{p}$ et sur les
fonctions polynômes définies dans $\mathfrak{p}$ joue un rôle
essentiel.

\subsection{}\label{1.8}
Dans \cite{parlett-2002}, Beresford Parlett démontre, pour l'essentiel
les formules \eqref{eq1.1} et \eqref{eq1.2} données plus haut :
\begin{align*}
  D(x)&=\det\Omega(x)\\
  D(x)&=\det g(x)
\end{align*}
qui expriment toutes deux le discriminant comme un déterminant, ce qui
justifie le titre de son article.

\subsection{}\label{1.9}
Dans l'article de Domokos cité (\cite{domokos-2011}), les méthodes (et
les résultats) franchissent un pas significatif vers le mieux possible
: Soit $V$ l'espace vectoriel introduit plus haut, i.e. celui
$\mathcal{J}(X)_{\frac{n(n-1)}{2}}$ des fonctions polynômes, homogènes
de degré $\frac{n(n-1)}{2}$, et nulles sur $X$. Cet espace joue déjà
un rôle dans l'article de Lax. Soit $W$ un sous-$SO(n)$-module de
$V$. Alors : quelle que soit la base orthonormée $(\varphi_{j})_{j}$
de $W$, la fonction $\sum_{j}(\varphi_{j})^{2}$ est
$SO(n)$-invariante, nulle sur $X$, et homogène de degré $n(n-1)$. Or
$D$ est le seul polynôme $SO(n)$-invariant, de degré $n(n-1)$, qui
s'annule sur $X$ (résultat déjà noté et utilisé par Lax). Donc :
$\sum_{j}(\varphi_{j})^{2}=cD$ (avec $c>0$).

Ainsi, chaque sous-$SO(n)$-module $W$ de $V$ fournit une décomposition de $D$ en
une somme de carrés, le nombre de carrés intervenant dans cette somme
étant la dimension de $W$. Comme Domokos poursuit le nombre minimum
$l(D)$ il aboutit au résultat suivant : $l(D)$ est majoré par le
minimum des dimensions des sous-espaces $SO(n)$-invariants de
$V$. Mieux : Domokos démontre que $V$ admet un sous-$SO(n)$-module
isomorphe à l'espace $H(n,n)$ des polynômes harmoniques de degré $n$
sur $\mathbb{R}^{n}$. D'où :
\begin{equation*}
l(D)\leq\dim H(n,n)
\end{equation*}
La détermination du minimum des dimensions des sous-espaces
$SO(n)$-invariants de $V$ reste à faire de même que celle de
$l(D)$. Dans le cas $n=3$, toutefois, Domokos montre que $l(D)$ vaut
$5$, alors que $\dim H(3,3)=7$.

\section{Ma contribution}\label{2}
Dans cette partie, je généralise ce qui précède dans le cadre des
espaces symètriques, mais la question plus subtile de la
détermination des nombres $l(D)$ pour les diverses fonctions $D$ ne
sera pas abordée. Pour cela, on utilisera librement les résultats
classiques de la théorie des algèbres de Lie semi-simples, qu'on
trouvera parmi d'autres dans les ouvrages de Helgason :
\cite{helgason-b-2008} et \cite{helgason-b-2001}.

Soit $\mathfrak{g}=\mathfrak{k}\oplus\mathfrak{p}$ une décomposition de
Cartan d'une algèbre de Lie {\color{black} semi-simple}
réelle $\mathfrak{g}$. On notera $G$ et $K$ des groupes de Lie
d'algèbres de Lie $\mathfrak{g}$ et $\mathfrak{k}$, et implicitement
$G$ et $K$ seront identifiés à des sous-groupes du groupe linéaire de
$GL(\mathfrak{g})$, via la représentation adjointe de $G$. Sur
$\mathfrak{g}$, il y a un produit scalaire $\langle\;\vert\;\rangle$
$K$-invariant qui a les propriétés suivantes :
\begin{equation*}
\langle x\vert y\rangle=-B(x,\theta y)
\end{equation*}
$\theta$ étant l'involution de Cartan associée à la décomposition du
même nom présentée plus haut, et $B$ étant la forme de Killing de
$\mathfrak{g}$. Pour ce produit scalaire, les espaces $\mathfrak{k}$
et $\mathfrak{p}$ sont orthogonaux, et la restriction de ce produit
scalaire à l'un et l'autre en fait des espaces euclidiens.

Comme dans l'introduction, on définit pour chaque élément $x$ de
$\mathfrak{p}$, les applications linéaires
$A(x):\mathfrak{p}\rightarrow\mathfrak{k}$ et
$B(x):\mathfrak{k}\rightarrow\mathfrak{p}$, chacune transformant un
élément $y$ en $[x,y]=\ad x.y$, et les composées $f(x)=B(x)\circ
A(x):\mathfrak{p}\rightarrow\mathfrak{p}$ et $g(x)=A(x)\circ
B(x):\mathfrak{k}\rightarrow\mathfrak{k}$, chacune transformant un
élément $y$ en $[x,[x,y]]=(\ad x)^{2}.y$. Comme $\ad x$ est une
transformation linéaire symétrique lorsque $x$ est un élément de
$\mathfrak{p}$ :
\begin{equation}
  \langle[x,\xi]\vert\eta\rangle=\langle\xi\vert[x,\eta]\rangle\mbox{ }
  \xi,\eta\in\mathfrak{g}
\end{equation}
il vient que $B(x)=A(x)^{*}$ est l'adjoint de $A(x)$, ou le transposé
de $A(x)$ lorsqu'on identifie $\mathfrak{k}$ et $\mathfrak{k}^{*}$
d'une part,  $\mathfrak{p}$ et $\mathfrak{p}^{*}$ d'autre part au
moyen du produit scalaire $\langle\;\vert\;\rangle$.

On réécrit les polynômes caractéristiques de $f(x)=B(x)\circ A(x)$ et
$g(x)=A(x)\circ B(x)$ :
\begin{align*}
  \det(tId_{\mathfrak{p}}-f(x))&=t^{r}-P_{1}(x)t^{r-1}+\cdots+(-1)^{r-l}P_{r-l}(x)
  t^{l}\\
  \det(tId_{\mathfrak{k}}-g(x))&=t^{s}-Q_{1}(x)t^{s-1}+\cdots+(-1)^{s-\alpha}Q_{s-\alpha}
  t^{\alpha}
\end{align*}
où $r=\dim\mathfrak{p}$, $s=\dim\mathfrak{k}$, $l$ est le rang de
$\mathfrak{p}$, et en vertu d'une propriété bien connue, on a :
\begin{equation*}
t^{s}\det(tId_{\mathfrak{p}}-f(x))=t^{r}\det(tId_{\mathfrak{p}}-g(x))
\end{equation*}
d'où : $\alpha=\dim\mathfrak{k}-\dim\mathfrak{p}+\rg\mathfrak{p}$
(=s-r+l).

\begin{pr}
Les $P_{j}$ ($1\leq j\leq r-l$) sont des sommes de carrés.
\end{pr}
\begin{dem}
Le produit scalaire $\langle\;\vert\;\rangle$ permet de définir un
produit scalaire naturel sur l'algèbre extérieure
$\wedge\mathfrak{g}$, en particulier sur chaque
$\wedge^{k}\mathfrak{g}$, les $\wedge^{k}\mathfrak{g}$ étant deux à
deux orthogonaux. Pour celà, on peut procéder comme suit :

Soit $(a_{i})_{1\leq i\leq r}$ une base orthonormée de
$\mathfrak{p}$, et soit $(b_{j})_{1\leq j\leq s}$ une base orthonormée de
$\mathfrak{k}$, de sorte que $(\alpha_{i})_{1\leq i\leq r+s}$, avec :
\begin{equation*}
\alpha_{i}=a_{i}\mbox{ }(1\leq i\leq r)\mbox{,
}\alpha_{r+j}=b_{j}\mbox{ }(1\leq j\leq s)
\end{equation*}
est une base orthonormée de $\mathfrak{g}$. Ceci dit
$(\alpha_{I})_{I}$, où $I=(i_{1},i_{2},\ldots,i_{k})$ avec : $1\leq
i_{1}<i_{2}<\cdots<i_{k}\leq r+s$ et
$\alpha_{I}=\alpha_{i_{1}}\wedge\alpha_{i_{2}}\wedge\cdots\wedge\alpha_{i_{k}}$,
est une base de $\wedge^{k}\mathfrak{g}$. On définit la structure
euclidienne sur $\wedge^{k}\mathfrak{g}$ en déclarant que
$(\alpha_{I})_{I}$ est une base orthonormée. En particulier,
$(a_{I})_{I}$ (où $a_{I}=a_{i_{1}}\wedge\cdots\wedge a_{i_{k}}$) est
une base orthonormée de $\wedge^{k}\mathfrak{p}$.

Il est bien connu que $P_{k}(x)$ est la trace de l'endomorphisme
$\wedge^{k}f(x)$ de $\wedge^{k}\mathfrak{p}$. Par suite :
\begin{equation*}
  P_{k}(x)=\sum_{I}\langle\wedge^{k}f(x)a_{I}\vert a_{I}\rangle
\end{equation*}
Mais $\wedge^{k}f(x)=\wedge^{k}B(x)\circ\wedge^{k}A(x)$ et
$B(x)=A(x)^{*}$, et vu la construction du produit scalaire sur
$\wedge^{k}\mathfrak{p}$ et $\wedge^{k}\mathfrak{k}$, on vérifie que
$\wedge^{k}A(x)^{*}=(\wedge^{k}A(x))^{*}$, et on arrive à ce qui nous
intéresse :
\begin{equation}
P_{k}(x)=\sum_{I}\Vert\wedge^{k}\ad x.a_{I}\Vert^{2}
\end{equation}
où il apparait que $P_{k}(x)$ est la somme des carrés des normes
euclidiennes des vecteurs $\wedge^{k}\ad x.a_{I}$, et chacun de ces
carrés est la somme de $
\begin{pmatrix}
  s\\
  k
\end{pmatrix}$ carrés de polynômes.
\end{dem}
\begin{rem}
D'après la démonstration de la proposition précédente, $P_{k}$ est
somme de $\begin{pmatrix} s\\ k
\end{pmatrix}\times\begin{pmatrix}
  r\\
  k
\end{pmatrix}$ carrés de polynômes, et ce nombre est scandaleusement
trop grand, au moins en comparaison des résultats obtenus par Domokos.
\end{rem}

\subsection{Un complément :}\label{2.1}
$f(x):\mathfrak{p}\mapsto\mathfrak{p}$ est un endomorphisme symétrique
de l'espace euclidien $\mathfrak{p}$, semi-défini positif. Pour de
tels endomorphismes, le rang se lie sur le polynôme caractéristique,
précisément :
\begin{equation*}
P_{k}(x)=0 \mbox{ si et seulement si le rang de $f(x)$ est strictement
inférieur à $k$}
\end{equation*}
De plus le rang de $f(x)=A(x)^{*}A(x)$ est celui de
$A(x)^{*}=B(x)$. Comme $B(x):\mathfrak{k}\rightarrow\mathfrak{p}$ est
l'application $z\mapsto[x,z]$, son image n'est autre que l'espace
tangent au point $x$ de l'orbite de $x$ sous
$K$. Conclusion :
\begin{equation*}
P_{k}(x)=0 \mbox{ si et seulement si $\dim K.x<k$}
\end{equation*}
On voit que les $P_{k}$ sont un cas particulier des polynômes sommes
de carrés dont l'expression est prévue par Domokos dans sa proposition
2.2.

Il y a en effet un procédé général pour construire des polynômes
invariants qui sont sommes de carrés.

\subsection{Covariants et sommes de carrés}\label{2.2}
Soient $V$ et $W$ deux $K$-modules unitaires, et soit
$\phi:W\rightarrow V$ un covariant, i.e. une fonction polynômiale
$K$-invariante au sens suivant :
\begin{equation*}
\phi(g.x)=g.\phi(x)\mbox{ }(x\in W\mbox{, }g\in K)
\end{equation*}
Alors la fonction $F:W\rightarrow\mathbb{R}$ définie ainsi :
\begin{equation*}
  F(x)=\Vert\phi(x)\Vert^{2}
  \mbox{ }(=\langle\phi(x)\vert\phi(x)\rangle)
\end{equation*}
est $K$-invariante, polynômiale, et somme des carrés des coordonnées
de $\phi(x)$ (relativement à une base orthonormée de $V$).

Un cas particulier, ou plutôt un exemple de cette construction,
redonne la proposition 2.2 de Domokos : Soit $\rho: K\rightarrow U(V)$
une représentation unitaire de $K$ (on notera encore $\rho$ la
représentation (dérivée) de $\mathfrak{k}$ dans $V$). Pour chaque $x$
dans $V$, on définit (l'application \flqq\;orbitale infinitésimale \frqq)
$\varphi_{x}:\mathfrak{k}\rightarrow V$ par :
\begin{equation*}
  \varphi_{x}(X)=\rho(X).x\mbox{ }(X\in\mathfrak{k})
\end{equation*}
La fonction $\phi:V\rightarrow \Hom_{\mathbb{R}}(\mathfrak{k},V)$
définie par : $\phi(x)=\varphi_{x}$ est un covariant, l'action de $K$
dans $\Hom_{\mathbb{R}}(\mathfrak{k},V)$ étant celle naturelle associée
à celles de $K$ dans $\mathfrak{k}$ (la représentation adjointe) et de
$K$ dans $V$ (la représentation $\rho$). Ceci étant, soit $m$ un
entier vérifiant : $1\leq m\leq\dim V$. On définit
$\phi_{m}:V\rightarrow \Hom_{\mathbb{R}}(\mathfrak{k}^{m},\wedge^{m}V)$
de la manière suivante :
\begin{equation*}
  \phi_{m}(x).(X_{1},\ldots,X_{m})=
  \varphi_{x}(X_{1})\wedge\cdots\wedge\varphi_{x}(X_{m})=
  \rho(X_{1}).x\wedge\cdots\wedge\rho(X_{m}).x
\end{equation*}
et on vérifie que $\phi_{m}$ est un covariant, relativement à la
représentation $\rho$ dans $V$ et à la représentation naturelle de $K$
dans $\mathfrak{k}^{m}\otimes\wedge^{m}V$. On voit que :
\begin{equation*}
  \phi_{m}(x)=0\mbox{ si et seulement si le rang de $\varphi_{x}$ est
    strictement inférieur à $m$}
\end{equation*}
et $F_{m}(x)=\Vert\phi_{m}(x)\Vert^{2}$ est un polynôme invariant tel
que :
\begin{equation*}
F_{m}(x)=0\mbox{ si et seulement si }\dim K.x<m 
\end{equation*}
Ces $F_{m}$ sont les polynômes invariants construits par Domokos.

Ces sommes de carrés font intervenir des mineurs de matrices. On
explique ceci dans le cas des polynômes $P_{k}$ introduits plus
haut. Pour celà, on introduit la matrices $M(x)$ de l'application
$A(x):\mathfrak{p}\rightarrow\mathfrak{k}$ relativement aux bases
$(a_{i})_{i}$ et $(b_{j})_{j}$ de $\mathfrak{p}$ et
$\mathfrak{k}$. Alors : $P_{k}(x)$ est la somme des carrés des
$k\times k$ mineurs de la matrice $M(x)$. Et clairement, on peut
remplacer $M(x)$ par sa transposée $M(x)^{*}$, qui est celle de
$B(x):\mathfrak{k}\rightarrow\mathfrak{p}$. On retrouve ainsi la
relation entre les $P_{k}$ et les dimensions des $K$-orbites.

\subsection{Des exemples}\label{2.3}


\subsubsection{La paire $(gl(n,\mathbb{R}),so(n))$}
On reprend brièvement l'exemple où se situent les travaux de
  Lax, Domokos et autres cités dans l'introduction. Ici :
\begin{equation*}
\mathfrak{g}=gl(n,\mathbb{R})=so(n)\oplus\mathfrak{p}
\end{equation*}
$s=\frac{n(n-1)}{2}$, $r=\frac{n(n+1)}{2}$ est la dimension de
l'espace vectoriel des matrices symétriques $\mathfrak{p}$. On note
$\mathfrak{a}$ le sous-espace de $\mathfrak{p}$ constitué par les
matrices diagonales. Dans la théorie générale des espaces symétriques,
$\mathfrak{a}$ est ce qu'on appelle un sous-espace de Cartan, et sa
dimension est le rang $l$ de la paire symétrique
$(\mathfrak{g},\mathfrak{p})$. Ici $l=n$. Une propriété importante des
sous-espaces de Cartan est qu'ils rencontrent toutes les $K$-orbites ;
il en résulte qu'une fonction $K$-invariante sur $\mathfrak{p}$ est
entièrement déterminée par sa restriction à un sous-espace de
Cartan. Dans cet exemple, on reconnait la propriété selon laquelle les
matrices symétriques sont diagonalisables par $SO(n)$.

On note $E'_{ij}$ la matrice (symétrique)
$\frac{1}{2}(E_{ij}+E_{ji})$, et $E''_{ij}$ (où $i<j$) la matrice
(antisymétrique) $\frac{1}{2}(E_{ij}-E_{ji})$. Avec :
\begin{equation*}
x=\diag(\lambda_{1},\ldots,\lambda_{n})
\end{equation*}
on a : $[x,E'_{ij}]=(\lambda_{i}-\lambda_{j})E''_{ij}$ (pour $i\leq
j$), $[x,E''_{ij}]=(\lambda_{i}-\lambda_{j})E'_{ij}$ (pour $i<j$) et :
\begin{align*}
f(x).E'_{ij}&=(\lambda_{i}-\lambda_{j})^{2}E'_{ij}\mbox{ }(i<j)\\
f(x).E_{ii}&=0\\
g(x).E''_{ij}&=(\lambda_{i}-\lambda_{j})^{2}E''_{ij}\mbox{ }(i<j)
\end{align*}
On peut comparer les polynômes caractéristiques de $f(x)$ et $g(x)$ :
\begin{align*}
\det(tId_{\mathfrak{p}}-f(x))&=t^{n}\prod_{i<j}(t-(\lambda_{i}-\lambda_{j})^{2})\\
\det(tId_{\mathfrak{k}}-g(x))&=\prod_{i<j}(t-(\lambda_{i}-\lambda_{j})^{2})
\end{align*}
Le discriminant $D(x)$ vaut
$\prod_{i<j}(\lambda_{i}-\lambda_{j})^{2}$ et $g(x)$ est
inversible lorsque $D(x)\neq0$ :
\begin{equation*}
D(x)=\det g(x)
\end{equation*}

\subsubsection{La paire $(gl(n,\mathbb{C}),u(n))$}
Dans cet exemple $\mathfrak{g}$ est $gl(n,\mathbb{C})$ considérée
comme algèbre de Lie réelle (on restreint le corps des scalaires), et
la décomposition de Cartan considérée est celle-ci :
\begin{equation*}
gl(n,\mathbb{C})=\mathfrak{k}\oplus\mathfrak{p}=u(n)\oplus\mathfrak{p}
\end{equation*}
où $\mathfrak{k}$ est l'algèbre de Lie des matrices anti-hermitiennes,
i.e. celle du groupe $K=U(n)$, et $\mathfrak{p}$ est l'espace des
matrices hermitiennes. Comme sous-espace de Cartan, on utilisera
l'espace $\mathfrak{a}$ des matrices diagonales (réelles
puisqu'hermitiennes). On a :
\begin{equation*}
\dim\mathfrak{k}=\dim\mathfrak{p}=n^{2}\mbox{, }l=n
\end{equation*}
Comme base de $\mathfrak{p}$, on utilisera  :
\begin{equation*}
E_{ii}\mbox{ }(1\leq i\leq n)\mbox{, }E'_{ij}\mbox{,
}\sqrt{-1}E''_{ij}\mbox{ }(i<j)
\end{equation*}
et comme base de $\mathfrak{k}$, on utilisera :
\begin{equation*}
\sqrt{-1}E_{ii}\mbox{ }(1\leq i\leq n)\mbox{, }E''_{ij}\mbox{,
}\sqrt{-1}E'_{ij}\mbox{ }(i<j)
\end{equation*}
avec $x=\diag(\lambda_{1},\ldots,\lambda_{n})$ dans $\mathfrak{a}$, on
a :
\begin{equation*}
[x,E_{ii}]=0\mbox{,
}[x,E'_{ij}]=(\lambda_{i}-\lambda_{j})E''_{ij}\mbox{, }
[x,E''_{ij}]=(\lambda_{i}-\lambda_{j})E'_{ij}
\end{equation*}
Il vient immédiatement que $f(x)$ et $g(x)$ ont le même polynôme
caractéristique :
\begin{equation*}
t^{n}\prod_{i<j}(t-(\lambda_{i}-\lambda_{j})^{2})^{2}
\end{equation*}
et $D(x)=\prod_{i<j}(\lambda_{i}-\lambda_{j})^{4}$.

Le fait que $f(x)$ et $g(x)$ aient le même polynôme caractéristique
provient du fait que $\mathfrak{k}=\sqrt{-1}\mathfrak{p}$
; en notant $\varphi:\mathfrak{p}\rightarrow\mathfrak{k}$ la bijection
$y\mapsto\sqrt{-1}y$ on voit que : $g(x)\circ\varphi=\varphi\circ
f(x)$ pour tout $x$ dans $\mathfrak{p}$, ce qui explique l'égalité des
polynômes caractéristiques. Ceci arrive chaque fois que
$\mathfrak{g}=\mathfrak{k}\oplus\sqrt{-1}\mathfrak{k}$ est la
complexifiée de l'algèbre de Lie compacte $\mathfrak{k}$.

\subsubsection{Les paires $(so(p,q),so(p)\times so(q))$}
On examine ici quelques exemples des algèbres de Lie $so(p,q)$. Par
définition, $so(p,q)$ est l'algèbre de Lie du sous-groupe de
$GL(n,\mathbb{R})$ qui laisse invariante la forme quadratique (non
dégénérée) sur $\mathbb{R}^{n}$ :
\begin{equation*}
  Q(\xi_{1},\ldots,\xi_{n})=\xi_{1}^{2}+\cdots+\xi_{p}^{2}-\xi_{p+1}^{2}-
  \cdots-\xi_{n}^{2}
\end{equation*}
(où $n=p+q$ et $p\geq q$). Une matrice $X$ de taille $n$ est un
élément de $so(p,q)$ si et seulement si :
\begin{equation*}
X=
\begin{pmatrix}
  A&B\\
  C&D
\end{pmatrix}\mbox{ avec $A\in so(p)$, $B^{*}=C$, $D\in so(q)$}
\end{equation*}
Une décomposition de Cartan de $\mathfrak{g}$ est :
$\mathfrak{g}=\mathfrak{k}\oplus\mathfrak{p}$, avec :
\begin{align*}
\mathfrak{k}&=\left\{
\begin{pmatrix}
  A&0\\
  0&D
\end{pmatrix}\mbox{ ; }A\in so(p)\mbox{, }D\in so(q)
\right\}\simeq so(p)\times so(q)\\
\mathfrak{p}&=\left\{
\begin{pmatrix}
  0&B\\
  B^{*}&0
\end{pmatrix}\mbox{ ; }B\in M_{p,q}(\mathbb{R})
\right\}\simeq M_{p,q}(\mathbb{R})
\end{align*}
Comme sous-espace de Cartan, on utilise l'ensemble des matrices
$\displaystyle
\begin{pmatrix}
  0&B\\
  B^{*}&0
\end{pmatrix}$
avec :
\begin{equation*}
B= 
\begin{pmatrix}
\lambda_{1}&\cdots&0\\
\vdots&\ddots&\vdots\\
0&\cdots&\lambda_{q}\\
0&\cdots&0\\
\vdots&&\vdots\\
0&\cdots&0
\end{pmatrix}\mbox{,  }(\lambda_{1},\ldots,\lambda_{q})\in\mathbb{R}^{q}
\end{equation*}
où il apparait que le rang $l$ de la paire symétrique
$(\mathfrak{g},\mathfrak{k})$ est $q(=\min(p,q))$.

\subsubsubsection{Le cas $q=1$}
On examine le cas $\mathfrak{g}=so(p,1)$ (avec $p\geq3$) (l'algèbre
de Lie du groupe de Lorentz).

Ici $\mathfrak{p}$ s'identifie, comme espace euclidien, à
$\mathbb{R}^{p}$ muni du produit scalaire habituel, et le commutateur
$[x,y]$ de deux vecteurs de $\mathfrak{p}$, calculé dans l'algèbre de
Lie $\mathfrak{g}$ est la matrice antisymétrique $x\wedge y
=(x_{i}y_{j}-x_{j}y_{i})_{1\leq i,j\leq p}$, élément de $so(p)$ comme
on pouvait s'y attendre ; autrement dit : $A(x).y=x\wedge y$, où on
voit que si $x\neq0$ :
\begin{equation*}
\ker A(x)=\mathbb{R}x\mbox{ et }\rg A(x)=p-1
\end{equation*}
Poursuivant le calcul, on trouve :
\begin{equation*}
B(x)\circ A(x).y=f(x).y=\Vert x\Vert^{2}y-(x.^{t}x).y
\end{equation*}
où $x.^{t}x$ est le produit de matrices $M_{p,1}(\mathbb{R})
\times
M_{1,p}(\mathbb{R})\rightarrow M_{p}(\mathbb{R})$. Le polynôme
caractéristique de $f(x)$ est donc :
\begin{equation*}
  \det(tId_{\mathfrak{p}}-f(x))=t(t-\Vert x\Vert^{2})^{p-1}
\end{equation*}
Les valeurs propres de $f(x)$ sont : $0$ avec la multiplicité $1$ et
$\Vert x\Vert^{2}$ avec la multiplicité $p-1$ (lorsque $x\neq 0$
évidemment). On trouve donc les expressions des polynômes $P_{k}$
($1\leq k\leq p-1$) :
\begin{equation*}
P_{1}(x)=(p-1)\Vert x\Vert^{2}, P_{2}(x)=
\begin{pmatrix}
  p-1\\
  2
\end{pmatrix}
\Vert x\Vert^{4},\ldots,D(x)=P_{p-1}(x)=\Vert x\Vert^{2(p-1)}
\end{equation*}
On vérifie que :
\begin{equation*}
  \det(tId_{\mathfrak{k}}-g(x))=t^{\frac{(p-1)(p-2)}{2}}
  (t-\Vert x\Vert^{2})^{p-1}
\end{equation*}
et on notera que le polynôme caractéristique de
$\ad x:\mathfrak{g}\rightarrow\mathfrak{g}$, $(x\in\mathfrak{p})$ est :
\begin{equation*}
  \det(tId_{\mathfrak{g}}-\ad x)=t^{m}(t-\Vert x\Vert)^{p-1}
  (t+\Vert x\Vert)^{p-1}
\end{equation*}
où $m=1+\dim so(p-1)$ est la dimension du centralisateur de
$\mathfrak{a}=\mathbb{R}x$ dans $\mathfrak{g}$.

\subsubsubsection{Le cas $q=p$}
Ici $\mathfrak{g}=so(p,p)=\mathfrak{k}\oplus\mathfrak{p}$, où
$\mathfrak{k}$ s'identifie au produit $so(p)\times so(p)$ et
$\mathfrak{p}$ s'identifie à l'espace $M_{p}(\mathbb{R})$ des matrices
carrées de taille $p$. Avec :
\begin{equation*}
x=
\begin{pmatrix}
  0&X\\
  ^{t}X&0
\end{pmatrix}\mbox{, }X=\diag(x_{1},\ldots,x_{p})
\end{equation*}
on trouve :
\begin{align*}
  \det(tId_{\mathfrak{p}}-f(x))&
  =t^{p}\prod_{i<j}(t-(x_{i}-x_{j})^{2})(t-(x_{i}+x_{j})^{2})\\
 \det(tId_{\mathfrak{k}}-g(x))&
 =\prod_{i<j}(t-(x_{i}-x_{j})^{2})(t-(x_{i}+x_{j})^{2})\\
 \det(tId_{\mathfrak{g}}-\ad x)&
 =t^{p}\prod_{i<j}(t^{2}
 -(x_{i}-x_{j})^{2})(t^{2}-(x_{i}+x_{j})^{2})\\
 D(x)&=\prod_{i<j}(x_{i}^{2}-x_{j}^{2})^{2}
 =\det g(x)
\end{align*}
Contrairement au cas $so(p,1)$ (avec $p\geq3$), $g(x)$ est
génériquement bijective et le discriminant $D$ n'est autre que la
restriction à $\mathfrak{p}$ du discriminant de l'algèbre de Lie
$\mathfrak{g}$ défini par :
\begin{equation*}
  \det(tId_{\mathfrak{g}}-\ad
  X)=t^{l'}(t^{m-l'}-a_{1}(X)t^{m-l'-1}+\cdots+(-1)^{m-l'}a_{m-l'}(X))
\end{equation*}
le polynôme $a_{m-l'}$ étant non nul, et $l'$ étant le rang de
l'algèbre de Lie $\mathfrak{g}$ ; le discriminant de $\mathfrak{g}$
est, par définition, le polynôme $a_{m-l'}$.

\subsection{Les espaces symétriques de rang maximal}\label{2.4}
Il est temps d'expliquer la différence entre $so(p,1)$ et
$so(p,p)$. Soit $\mathfrak{g}=\mathfrak{k}\oplus\mathfrak{p}$ une
décomposition de Cartan de $\mathfrak{g}$ (l'application linéaire
$\theta:\mathfrak{g}\rightarrow\mathfrak{g}$ telle que
$\theta\vert_{\mathfrak{k}}=Id_{\mathfrak{k}}$ et
$\theta\vert_{\mathfrak{p}}=-Id_{\mathfrak{p}}$ est un automorphisme
d'ordre $2$ de l'algèbre de Lie $\mathfrak{g}$, une involution de
Cartan) et soit $\mathfrak{a}$ un sous-espace de Cartan de
$\mathfrak{p}$ ; par définition $\mathfrak{a}$ est un sous-espace
vectoriel de $\mathfrak{p}$, abélien au sens suivant : $[x,y]=0$ pour
tout couple $(x,y)$ d'éléments de $\mathfrak{a}$, et maximal par
rapport à cette propriété. La dimension $l$ de l'espace vectoriel
$\mathfrak{a}$ est le rang de l'espace symétrique $G/K$ ou de la paire
symétrique $(\mathfrak{g},\mathfrak{k})$, et il est clair que $l$ est
majoré par le rang de $\mathfrak{g}$ (qui est la dimension des
sous-algèbres de Cartan de $\mathfrak{g}$). Ceci explique qu'on dise
d'une paire symétrique $(\mathfrak{g},\mathfrak{k})$ qu'elle est de
rang maximal lorsque son rang coïncide avec celui de $\mathfrak{g}$,
i.e. lorsque un (ou tout) sous-espace de Cartan $\mathfrak{a}$ de
$\mathfrak{p}$ est une sous-algèbre de Cartan de $\mathfrak{g}$.

Soit $\mathfrak{a}$ un sous-espace de Cartan de $\mathfrak{p}$. Alors
$\mathfrak{a}$ est une sous-algèbre de Cartan de $\mathfrak{g}$ si et
seulement si $B(x):\mathfrak{k}\rightarrow\mathfrak{p}$ est
génériquement injective quand $x$ parcourt $\mathfrak{a}$, si et
seulement si $\det g(x)$ est un polynôme non nul.

Dans les exemples présentés ci-dessus, il y a deux paires symétriques
de rang maximal : $(gl(n,\mathbb{R}),so(n,\mathbb{R}))$ et
$(so(p,p),so(p)\times so(p))$, et deux paires de rang non maximal :
$(gl(n,\mathbb{C}),$ $u(n))$ et $(so(p,1),so(p))$. Parmi les algèbres
$so(p,q)$, seules $so(p,p)$ (cas où $q=p$) et $so(p+1,p)$ (cas où
$q=p-1$) donnent des paires symétriques de rang maximal.

\subsection{Cas général}\label{2.5}
On revient au cas général d'une paire symétrique
$(\mathfrak{g},\mathfrak{k})$, avec la décomposition de Cartan
$\mathfrak{g}=\mathfrak{k}\oplus\mathfrak{p}$. Soit $x$ dans
$\mathfrak{p}$, alors $\ad x$ (resp. $f(x)$, $g(x)$) est un
endomorphisme symétrique de $\mathfrak{g}$ (resp. $\mathfrak{p}$,
$\mathfrak{k}$). Soit $\lambda$ une valeur propre (évidemment réelle)
de $\ad x$.

\begin{itemize}[label=-]
\item Cas $\lambda=0$. L'espace propre associé à $\lambda$,
  i.e. $\ker\ad x$, est somme directe :
\begin{equation*}
\ker\ad x=(\ker\ad x)\cap\mathfrak{k}\oplus(\ker\ad x)\cap\mathfrak{p}
\end{equation*}

\item Cas $\lambda\neq0$. Lorsque $z+y$, avec $z\in\mathfrak{k}$ et
  $y\in\mathfrak{p}$ est un vecteur propre de $\ad x$ associé à
  $\lambda$, on a $[x,z]=\lambda y$, $[x,y]=\lambda z$,
  $f(x).y=\lambda^{2}y$, $g(x).z=\lambda^{2}z$ et $\lambda^{2}$ est
  une valeur propre de $f(x)$ (resp. $g(x)$). En sens inverse, lorsque
  $f(x).y=\mu y$ ($\mu>0$, $y\neq0$), alors $z+y$ où
  $z=\pm\frac{1}{\sqrt\mu}[x,y]$, est vecteur propre de $\ad x$,
  associé à la valeur propre $\pm\sqrt\mu$.
\end{itemize}

Soit $\mathfrak{a}$ un sous-espace de Cartan de $\mathfrak{p}$. Pour
chaque forme linéaire $\lambda$ sur $\mathfrak{a}$, non nulle, on note
$\mathfrak{g}_{\lambda}$ l'espace des $\xi$ dans $\mathfrak{g}$ tels
que $[x,\xi]=\lambda\xi$ pour tout $x$ dans $\mathfrak{a}$ (les $\ad
x$, $x$ parcourant $\mathfrak{a}$, sont simultanément
diagonalisables), et on note $\Delta(\mathfrak{a})$ l'ensemble des
$\lambda$ ci-dessus tels que $\mathfrak{g}_{\lambda}\neq\{0\}$ (c'est
l'ensemble des racines de $(\mathfrak{g},\mathfrak{a})$). Pour
$\lambda$ dans $\Delta(\mathfrak{a})$, on désigne par
$\mathfrak{p}_{\lambda}$ (resp. $\mathfrak{k}_{\lambda}$) l'ensemble
des $\xi$ dans $\mathfrak{p}$ (resp. $\mathfrak{k}$) tels que $(\ad
x)^{2}.\xi=\lambda^{2}\xi$ ; comme
$-\lambda=\theta(\lambda)=\lambda\circ\theta$ est une racine, on voit
que $\mathfrak{p}_{\lambda}=\mathfrak{p}_{-\lambda}$
(resp. $\mathfrak{k}_{\lambda}=\mathfrak{k}_{-\lambda}$). Ceci amène à
partager $\Delta(\mathfrak{a})$ en racines positives (dont l'ensemble
est noté $\Delta_{+}(\mathfrak{a})$) et racines négatives (constituant
$\Delta_{-}(\mathfrak{a})$) :
$\Delta(\mathfrak{a})=\Delta_{+}(\mathfrak{a})\sqcup\Delta_{-}(\mathfrak{a})$,
de sorte que des deux racines $\lambda$ et $-\lambda$, l'une est
positive et l'autre est négative. On vérifie que l'on a :
\begin{equation*}
\mathfrak{k}=
\mathfrak{m}\oplus\sum_{\lambda\in\Delta_{+}(\mathfrak{a})}\mathfrak{k}_{\lambda}
\mbox{, }
\mathfrak{p}=\mathfrak{a}\oplus
\sum_{\lambda\in\Delta_{+}(\mathfrak{a})}\mathfrak{p}_{\lambda}
\end{equation*}
où $\mathfrak{m}$ est le centralisateur de $\mathfrak{a}$ dans
$\mathfrak{k}$, et aussi :
\begin{equation*}
  \mathfrak{g}=
  \mathfrak{m}\oplus\mathfrak{a}\oplus\sum_{\lambda\in\Delta(\mathfrak{a})}
  \mathfrak{g}_{\lambda}
\end{equation*}
avec : $\mathfrak{g}_{\lambda}+\mathfrak{g}_{-\lambda}=
\mathfrak{p}_{\lambda}+\mathfrak{k}_{\lambda}$
($\lambda\in\Delta_{+}(\mathfrak{a})$).
On a donc (tout ceci avec $x$ dans $\mathfrak{a}$) :
\begin{equation*}
  \det(tId_{\mathfrak{g}}-\ad x)=t^{\dim\mathfrak{m}+l}
  \prod_{\lambda\in\Delta(\mathfrak{a})}(t-\lambda(x))^{\mu(\lambda)}
\end{equation*}
où
$\mu(\lambda)=\dim\mathfrak{g}_{\lambda}=\dim\mathfrak{p}_{\lambda}$
est la multiplicité de la racine $\lambda$.
\begin{align*}
  \det(tId_{\mathfrak{p}}-f(x))&=
  t^{l}\prod_{\lambda\in\Delta_{+}
    (\mathfrak{a})}(t-\lambda(x)^{2})^{\mu(\lambda)}\\
 \det(tId_{\mathfrak{k}}-g(x))&=
 t^{\dim\mathfrak{m}}\prod_{\lambda\in\Delta_{+}
   (\mathfrak{a})}(t-\lambda(x)^{2})^{\mu(\lambda)}\\
 D(x)&=\prod_{\lambda\in\Delta_{+}(\mathfrak{a})}
 \lambda(x)^{2\mu(\lambda)}
\end{align*}

Par définition, un élément $x$ de $\mathfrak{p}$ est dit régulier
lorsque $(\ker\ad x)\cap\mathfrak{p}$ est de dimension $l$, l'ensemble
des éléments réguliers est désigné par la notation
$\mathfrak{p}_{reg}$ ; l'ensemble
$\mathfrak{a}_{reg}=\mathfrak{p}_{reg}\cap\mathfrak{a}$ est celui des
$x$ dans $\mathfrak{a}$ tels que : $D(x)\neq0$, i.e. tels que :
\begin{equation*}
\lambda(x)\neq0\mbox{ pour toute racine $\lambda$}
\end{equation*}

\subsection{Exemples}\label{2.6} 

Dans ce numéro, nous appliquons les résultats de \ref{2.5} aux
exemples du numéro \ref{2.3}.

\subsubsection{La paire $(gl(n,\mathbb{R}),so(n))$}
L'espace $\mathfrak{p}$ est celui des matrices réelles symétriques,
$\mathfrak{a}$ est celui des matrices diagonales et
$\Delta(\mathfrak{a})$ est l'ensemble des $\lambda_{ij}$ ($i\neq j$) :
$\langle\lambda_{ij},\diag(x_{1},\ldots,x_{n})\rangle=x_{i}-x_{j}$,
$\mathfrak{g}_{\lambda_{ij}}=\mathbb{R}E_{ij}$,
$\mathfrak{p}_{\lambda_{ij}}=\mathbb{R}(E_{ij}+E_{ji})$,
$\mathfrak{k}_{\lambda_{ij}}=\mathbb{R}(E_{ij}-E_{ji})$,
$\mathfrak{m}=\{0\}$ et $\mu(\lambda)=1$ pour toute racine $\lambda$.
\begin{align*}
\det(tId_{\mathfrak{g}}-\ad x)&=t^{n}\prod_{i\neq
  j}(t-(x_{i}-x_{j}))\\
\det(tId_{\mathfrak{p}}-f(x))&=t^{n}\prod_{i<
  j}(t-(x_{i}-x_{j})^{2})
\end{align*}
$\Delta_{+}(\mathfrak{a})$ est l'ensemble des $\lambda_{ij}$, $i<j$,
et :
\begin{equation*}
D(x_{1},\ldots,x_{n})=\prod_{i<j}(x_{i}-x_{j})^{2}
\end{equation*}

\subsubsection{La paire $(gl(n,\mathbb{C}),u(n))$}
L'espace $\mathfrak{p}$ est celui des matrices hermitiennes,
$\mathfrak{a}$ est celui des matrices diagonales réelles et :
$\Delta(\mathfrak{a})$ est l'ensemble des $\lambda_{ij}$ ($i\neq j$),
$\langle\lambda_{ij},\diag(x_{1},\ldots,x_{n})\rangle=x_{i}-x_{j}$,
\begin{align*}
  \mathfrak{g}_{\lambda_{ij}}&=\mathbb{R}E_{ij}\oplus\mathbb{R}\sqrt{-1}
  E_{ij}\mbox{ }(\mu(\lambda_{ij})=2)
  \\
 \mathfrak{p}_{\lambda_{ij}}&=\mathbb{R}{E}'_{ij}\oplus\mathbb{R}\sqrt{-1}
  E''_{ij}\mbox{ }(i<j)
\end{align*}
$\mathfrak{m}$ est le sous-espace des matrices diagonales
antihermitiennes
\begin{equation*}
  D(x)=\prod_{i<j}(x_{i}-x_{j})^{4}
\end{equation*}
Noter qu'ici : $\mathfrak{m}\oplus\mathfrak{a}$ est une sous-algèbre
de Cartan de $\mathfrak{g}=gl(n,\mathbb{C})$, dont le rang est $2n$,
tandis que celui de l'espace symétrique $GL(n,\mathbb{C})/U(n)$ est $n$.

\subsubsection{Paires orthogonales}

\subsubsubsection{La paire $(so(p,1),so(p))$}
L'espace $\mathfrak{p}$ a été identifié plus haut à $\mathbb{R}^{p}$,
$\mathfrak{a}$ est la droite $\mathbb{R}e_{1}$
($e_{1}=(1,0,\ldots,0)$) et : $\Delta(\mathfrak{a})=\{\pm\lambda\}$
avec $\langle\lambda,(\xi,0,\ldots,0)\rangle=\xi$.

Soit $x$ dans $\mathfrak{a}$, avec $x=\xi e_{1}$, $\xi\neq0$, de sorte
que $\lambda(x)=\xi$. Pour définir $\mathfrak{g}_{\xi}$ et
$\mathfrak{g}_{-\xi}$, on définit
$\varphi_{\pm}:(\mathbb{R}x)^{\perp}\rightarrow\mathfrak{g}$ :
$\varphi_{\pm}(y)=\pm\frac{1}{\xi}x\wedge y+y$ et alors :
\begin{equation*}
\mathfrak{g}_{\xi}=\im\varphi_{+}\mbox{, }\mathfrak{g}_{-\xi}=\im\varphi_{-}
\end{equation*}
Donc $\mu(\xi)=p-1=\mu(-\xi)$. De plus :
\begin{equation*}
\mathfrak{p}_{\xi}=(\mathbb{R}x)^{\perp}
\end{equation*}
Ici $\mathfrak{m}$ est la sous-algèbre de Lie de $so(p)$, qui
\flqq\;stabilise \frqq~ $x$, i.e. constituée par les matrices $Z$ dans
$so(p)$ telles que $Z.x=0$. 
\begin{equation*}
D(x)=\Vert x\Vert^{2(p-1)}
\end{equation*}

En résumé avec $m'=1+\dim so(p-1)$ :
\begin{align*}
  \det(tId_{\mathfrak{g}}-\ad
  x)&=t^{m'}(t-\lambda(x))^{p-1}(t+\lambda(x))^{p-1}\\
  \det(tId_{\mathfrak{p}}-f(x))&=t(t-\lambda(x)^{2})^{p-1}\\
  \det(tId_{\mathfrak{k}}-g(x))&=t^{\dim so(p-1)}
  (t-\lambda(x)^{2})^{p-1} 
\end{align*}

\subsubsubsection{La paire $(so(p,p),so(p)\times so(p))$} On a
identifié $\mathfrak{p}$ à l'espace $M_{p}(\mathbb{R})$, et
$\mathfrak{a}$ à celui des matrices diagonales de
$M_{p}(\mathbb{R})$. 
\begin{equation*}
  \Delta(\mathfrak{a})=\{\pm\lambda_{ij},\pm\lambda'_{ij}\mbox{ ;
  }i<j\}
\end{equation*}
  \begin{align*}
  \langle\lambda_{ij},\diag(x_{1},\ldots,x_{p})\rangle&=x_{i}-x_{j}\\
\langle\lambda'_{ij},\diag(x_{1},\ldots,x_{p})\rangle&=x_{i}+x_{j}  
\end{align*}
Avec, comme plus haut :
\begin{equation*}
2E'_{ij}=E_{ij}+E_{ji}\mbox{, }2E''_{ij}=E_{ij}-E_{ji}\mbox{ }(1\leq
i<j\leq p)
\end{equation*}
on a :
\begin{alignat*}{2}
  \mathfrak{g}_{\lambda_{ij}}&=\mathbb{R}
\begin{pmatrix}
  E''_{ij}&E'_{ij}\\
  E'_{ij}&E''_{ij}
\end{pmatrix} & \qquad
\mathfrak{g}_{-\lambda_{ij}}&=\mathbb{R}
\begin{pmatrix}
  E''_{ij}&-E'_{ij}\\
  -E'_{ij}&E''_{ij}
\end{pmatrix} \\
\mathfrak{g}_{\lambda'_{ij}}&=\mathbb{R}
\begin{pmatrix}
  E''_{ij}&E'_{ij}\\
  E'_{ij}&-E''_{ij}
\end{pmatrix} & \qquad
\mathfrak{g}_{-\lambda'_{ij}}&=\mathbb{R}
\begin{pmatrix}
  E''_{ij}&-E'_{ij}\\
  -E'_{ij}&-E''_{ij}
\end{pmatrix} 
\end{alignat*}
et on voit que $\mu(\lambda)=1$ pour tout $\lambda$ dans
$\Delta(\mathfrak{a})$, et :
\begin{equation*}
  \det(tId_{\mathfrak{g}}-\ad x)=t^{p}\prod_{i<j} (t^{2}-\lambda_{ij}(x)^{2})
(t^{2}-\lambda'_{ij}(x)^{2})  
\end{equation*}
Pour les valeurs propres de $f(x)$ :
\begin{equation*}
\mathfrak{p}_{\lambda_{ij}}=\mathbb{R}
\begin{pmatrix}
  0&E'_{ij}\\
  E'_{ij}&0
\end{pmatrix} \qquad
\mathfrak{p}_{\lambda'_{ij}}=\mathbb{R}
\begin{pmatrix}
  0&E''_{ij}\\
  -E''_{ij}&0
\end{pmatrix}
\end{equation*}
\begin{equation*}
\det(tId_{\mathfrak{p}}- f(x))=t^{p}\prod_{i<j} (t-\lambda_{ij}(x)^{2})
(t-\lambda'_{ij}(x)^{2})
\end{equation*}
Enfin ($\Delta_{+}$ est l'ensemble des $\lambda_{ij}$, $\lambda'_{ij}$
avec $i<j$)
\begin{equation*}
\det(tId_{\mathfrak{k}}-g(x))=\prod_{i<j} (t-\lambda_{ij}(x)^{2})
(t-\lambda'_{ij}(x)^{2})
\end{equation*}

\subsection{L'héritage d'Ilyushechkin}
Notons $\mathbb{R}[\mathfrak{p}]^{K}$ l'algèbre des fonctions
polynômes (sur $\mathfrak{p}$) qui sont $K$-invariantes. Il est connu
que $\mathbb{R}[\mathfrak{p}]^{K}$ est une algèbre (pure) de polynômes
de dimension (de Krull) $l$ (i.e. de dimension égale au rang de la
paire $(\mathfrak{g},\mathfrak{k})$, et aussi à la dimension de
l'espace vectoriel $\mathfrak{a}$).

Soit $(p_{1},\ldots,p_{l})$ un système de générateurs (homogènes)
agébriquement indépendants de l'algèbre
$\mathbb{R}[\mathfrak{p}]^{K}$. On supposera que les $p_{j}$ sont
numérotés de telle façon que leurs degrés sont rangés dans l'ordre
croissant : $m_{1}+1=\deg p_{1}\leq m_{2}+1=\deg p_{2}\leq\cdots\leq
m_{l}+1=\deg p_{l}$ ; les $m_{j}$ sont ce qu'on appelle les exposants
de la paire $(\mathfrak{g},\mathfrak{k})$.

Pour chaque fonction $f:\mathfrak{p}\rightarrow\mathbb{R}$, on notera
$\nabla f:\mathfrak{p}\rightarrow\mathfrak{p}$ le gradient de $f$
calculé au moyen du produit scalaire $\langle\;\vert\;\rangle$ :
\begin{equation*}
  \left(\frac{\diff}{\diff t}\right)_{0}f(x+ty)=\langle\diff f(x),y\rangle=
  \langle\nabla f(x)\vert y\rangle
\end{equation*}

Soit $F:\mathfrak{p}\rightarrow\wedge^{l}\mathfrak{p}$ définie par : 
\begin{equation*}
F(x)=\nabla p_{1}(x)\wedge\cdots\wedge\nabla p_{l}(x)\mbox{ }(x\in\mathfrak{p})
\end{equation*}
C'est une fonction polynômiale, homogène de degré :
$m_{1}+m_{2}+\cdots+m_{l}$, et comme chaque $\nabla p_{j}$ est un
covariant ($\nabla p_{j}(g.x)=g.\nabla p_{j}(x)$ pour tout $g$ dans
$K$) il vient que $F$ est un covariant :
\begin{equation*}
F(g.x)=(\wedge^{l}g).F(x)
\end{equation*}
et comme expliqué plus haut, la fonction :
\begin{equation*}
\phi(x)=\Vert F(x)\Vert^{2}
\end{equation*}
définit sur $\mathfrak{p}$ une fonction polynôme $K$-invariante,
homogène de degré $2(m_{1}+\cdots+m_{l})$, qui est une somme de
carrés.

Il y a le théorème de Chevalley : L'opération $Res_{\mathfrak{a}}$ qui
restreint à $\mathfrak{a}$ les fonctions définies sur $\mathfrak{p}$
est un isomorphisme d'algèbres de $\mathbb{R}[\mathfrak{p}]^{K}$ sur
l'algèbre $\mathbb{R}[\mathfrak{a}]^{W}$ des fonctions polynômes
$W$-invariantes sur $\mathfrak{a}$. Le groupe (de Weyl) $W$ se définit
comme suit : On introduit le sous-groupe de $K$ qui \flqq\;
normalise \frqq~ $\mathfrak{a}$ :
\begin{equation*}
  Norm_{K}(\mathfrak{a})=\{g\in K\;\vert\;g.\mathfrak{a}\subset\mathfrak{a}\}
\end{equation*}
Puis, chaque $g$ dans $Norm_{K}(\mathfrak{a})$ définit une
transformation linéaire bijective de $\mathfrak{a}$ : $x\mapsto g.x$
($x\in\mathfrak{a}$), et $W$ est le sous-groupe de $GL(\mathfrak{a})$
constitué par toutes les transformations linéaires ainsi
obtenues. Ceci étant, on note $q_{j}$ la restriction à $\mathfrak{a}$
de $p_{j}$, de sorte que 
\begin{equation*}
\mathbb{R}[\mathfrak{a}]^{W}=\mathbb{R}[q_{1},\ldots,q_{l}]
\end{equation*}

Chaque $q_{j}$ admet un gradient $\nabla
q_{j}:\mathfrak{a}\rightarrow\mathfrak{a}$ que l'on calcule au moyen
du produit scalaire de $\mathfrak{a}$, et il se trouve que lorsque $x$
est dans $\mathfrak{a}$, $\nabla p_{j}(x)$ appartient à $\mathfrak{a}$
et coïncide avec le vecteur $\nabla q_{j}(x)$ :
\begin{equation*}
  \nabla q_{j}(x)=\nabla p_{j}(x)\mbox{ pour tout $x$ dans $\mathfrak{a}$}
\end{equation*}
Par suite, lorsque $x\in\mathfrak{a}$ :
\begin{equation*}
F(x)=\nabla p_{1}(x)\wedge\cdots\wedge\nabla p_{l}(x)
\end{equation*}
appartient au sous-espace $\wedge^{l}\mathfrak{a}$ de
$\wedge^{l}\mathfrak{p}$ et :
\begin{align*}
  F(x)&=\nabla q_{1}(x)\wedge\cdots\wedge\nabla q_{l}(x)\\
  \Vert F(x)\Vert^{2}&=\det\Omega(x)
\end{align*}
avec $\Omega(x)=(\langle\nabla q_{i}(x)\vert\nabla
q_{j}(x)\rangle)_{1\leq i,j\leq l}$. En effet, si on fixe une base
orthonormale $(e_{1},\ldots,e_{l})$ de $\mathfrak{a}$ et si on designe
par $\partial_{i}$ la dérivée le long du vecteur $e_{i}$, alors :
\begin{equation*}
  \nabla q_{1}(x)\wedge\cdots\wedge\nabla q_{l}(x)=
  \det(Jac(q_{1},\ldots,q_{l}))\omega
\end{equation*}
où $\omega=e_{1}\wedge\cdots\wedge e_{l}$ et $Jac(q_{1},\ldots,q_{l})$
est la matrices jacobienne $(\partial_{i}q_{j})_{ij}$, évaluée au point $x$.
Vu la $K$-invariance de $\phi$, on conclut :
\begin{equation*}
\phi(x)=\det((\langle\nabla p_{i}(x)\vert\nabla
p_{j}(x)\rangle)_{1\leq i,j\leq l})
\end{equation*}
pour tout $x\in\mathfrak{p}$.

\subsection{Compléments} Le groupe $W$  est un sous-groupe fini de
$GL(\mathfrak{a})$, engendré par des réflections. L'ordre de $W$ et le
nombre de réflections de $W$ se lisent sur les degrés
$d_{1},\ldots,d_{l}$ des invariants \flqq\;fondamentaux \frqq~
$p_{1},\ldots,p_{l}$ :

(i) l'ordre de $W$ est le produit $d_{1}d_{2}\cdots d_{l}$

(ii) le nombre de réflections est $\sum_{j}(d_{j}-1)$

Lorsque $x$ est dans $\mathfrak{a}$, en notant $F_{0}$ la restriction
de $F$ à $\mathfrak{a}$, on a
\begin{equation*}
F_{0}(x)=Jac(q_{1},\ldots,q_{l})e_{1}\wedge\cdots\wedge e_{l}
\end{equation*}
On simplifie les notations en écrivant :
\begin{equation*}
J(x)=Jac(q_{1},\ldots,q_{l})(x)
\end{equation*}
Alors $J$ est un polynôme semi-invariant (sous l'action de $W$) :
\begin{equation*}
J(w.x)=(\det w)J(x)\mbox{ }(w\in W\mbox{, }x\in\mathfrak{a})
\end{equation*}
(et : $\phi(x)=J(x)^{2}$), et son degré est le nombre de réflections
présentes dans $W$.

Par ailleurs $F(x)=\diff p_{1}(x)\wedge\cdots\wedge\diff p_{l}(x)$ est
non nul si et seulement si les vecteurs $\diff p_{l}(x),\diff
p_{2}(x),\ldots,\diff p_{l}(x)$ sont linéairement indépendants ; de même,
$F_{0}(x)\neq0$ si et seulement si $\diff q_{1}(x),\diff
q_{2}(x),\ldots,\diff q_{l}(x)$ forment une base de $\mathfrak{a}$, si
et seulement si $x$ est un élément régulier. Donc $J(x)=0$ lorsque
$\lambda(x)=0$, $\lambda$ étant une racine dans
$\Delta_{+}(\mathfrak{a)}$. Ainsi, $J$ est divisible par toute racine
positive.

\subsection{Exemples}

\subsubsection{ Le cas $(gl(n,\mathbb{R}),so(n))$}
$\mathbb{R}[\mathfrak{p}]^{K}=\mathbb{R}[p_{1},p_{2},\ldots,p_{n}]$, 
avec : $p_{j}(x)=\frac{1}{j}\Tr x^{j}$ pour tout $x$ dans
$\mathfrak{p}$. Les restrictions $q_{j}$ des $p_{j}$ à $\mathfrak{a}$
sont définies ainsi :
\begin{equation*}
  q_{j}(\diag(x_{1},\ldots,x_{n}))=
  q_{j}(x_{1},\ldots,x_{n})=\frac{1}{j}\sum_{i}x_{i}^{j}
\end{equation*}
et $\nabla q_{j+1}(x_{1},\ldots,x_{n})=(x_{1}^{j},x_{2}^{j},
,\ldots,x_{n}^{j})$. Donc $J(x)$ est le discriminant de la matrice de
Van der Monde : 
\begin{equation*}
\begin{vmatrix}
  1&x_{1}&\cdots&x_{1}^{n-1}\\
  1&x_{2}&\cdots&x_{2}^{n-1}\\
  \vdots&\vdots& &\vdots\\
  1&x_{n}&\cdots&x_{n}^{n-1}
\end{vmatrix}
\end{equation*}
et :
$J(x)=\prod_{i<j}(x_{i}-x_{j})=\prod_{\lambda\in\Delta_{+}(\mathfrak{a})}\lambda(x)$.
Dans cet exemple :
\begin{equation*}
J^{2}(x)=\prod_{\lambda\in\Delta_{+}(\mathfrak{a})}\lambda^{2}(x)=D(x)\mbox{
  }(x\in\mathfrak{a})
\end{equation*}
et par suite, pour tout $x\in\mathfrak{p}$ :
\begin{equation*}
D(x)=\Vert\nabla p_{1}(x)\wedge\cdots\wedge\nabla p_{n}(x)\Vert^{2}
\end{equation*}

\subsubsection{Le cas $(so(p,1),so(p))$}
$\mathbb{R}[\mathfrak{p}]^{K}=\mathbb{R}[Q]$,
avec : $Q(x_{1},\ldots,x_{p})=\frac{1}{2}(x_{1}^{2}+\cdots+x_{p}^{2})$,
sachant que $\mathfrak{p}\simeq\mathbb{R}^{p}$. La restriction  $q$
à $\mathfrak{a}=\mathbb{R}e_{1}$ est donc :
\begin{equation*}
  q(x_{1})=Q(x_{1},0,\ldots,0))=\frac{1}{2}x^{2}
\end{equation*}
\begin{equation*}
  \nabla q(x_{1})=x_{1}\mbox{, }J(x)=x_{1}
\end{equation*}
\begin{equation*}
J^{2}(x)=\Vert x\Vert^{2}
\end{equation*}
Dans cet exemple, lorsque $x\in\mathfrak{a}$ :
\begin{equation*}
  J^{2}(x)=\prod_{\lambda\in\Delta_{+}(\mathfrak{a})}\lambda^{2}(x)
\end{equation*}
mais $\Vert\nabla Q(x)\Vert^{2}=\Vert x\Vert^{2}$ tandis que :
\begin{equation*}
D(x)=\Vert x\Vert^{2(p-1)}
\end{equation*}

\subsubsection{Le cas $(so(p,p),so(p)\times so(p))$ }
$\mathbb{R}[\mathfrak{p}]^{K}=\mathbb{R}[P_{1},\ldots,P_{p}]$,
avec : $P_{j}(x)=\frac{1}{2j}\Tr((xx^{*})^{j})$, $1\leq j\leq p-1$,
$P_{p}(x)=\det x$, de sorte que les restrictions $Q_{j}$ des $P_{j}$
aux matrices diagonales sont données par :
\begin{align*}
Q_{j}(\diag(x_{1},\ldots,x_{p}))&=\frac{1}{2j}(x_{1}^{2j}+\cdots+x_{p}^{2j})\mbox{
}(1\leq j\leq p-1)\\
Q_{p}(\diag(x_{1},\ldots,x_{p}))&=x_{1}x_{2}\cdots x_{p}
\end{align*}
Dès lors, la matrice jacobienne des $Q_{j}$ se présente sous la forme
:
\begin{equation*}
\begin{vmatrix}
  x_{1}&x_{2}&\cdots&x_{p}\\
  x_{1}^{3}&x_{2}^{3}&\cdots&x_{p}^{3}\\
  \vdots&\vdots& &\vdots\\
  x_{1}^{2p-3}&x_{2}^{2p-3}&\cdots&x_{p}^{2p-3}\\
  x_{2}x_{3}\cdots x_{p}&x_{1}x_{3}\cdots
  x_{p}&\cdots&x_{1}x_{2}\cdots x_{p-1}
\end{vmatrix}
\end{equation*}
et un calcul élémentaire donne, au signe près :
\begin{equation*}
J(x_{1},\ldots,x_{p})=\prod_{i<j}(x_{i}^{2}-x_{j}^{2})
\end{equation*}
Dans cet exemple :
\begin{equation*}
D(x)=\Vert\nabla P_{1}(x)\wedge\cdots\wedge\nabla P_{p}(x)\Vert^{2}
\end{equation*}
pour tout $x$ dans $\mathfrak{p}$ (à une multiplication près par un
nombre réel strictement positif).

\subsubsection{Le cas $(gl(n,\mathbb{C}),u(n))$}
$\mathbb{R}[\mathfrak{p}]^{K}=\mathbb{R}[p_{1},p_{2},\ldots,p_{n}]$, 
avec : $p_{j}(x)=\frac{1}{j}\Tr x^{j}$ pour toute matrice hermitienne
$x$ (ce qui explique que les $p_{j}$ soient à valeurs réelles). Les
restrictions $q_{j}$ aux matrices diagonales sont celles trouvées dans
l'exemple $(gl(n,\mathbb{R}),so(n))$ et par suite :
\begin{gather*}
  J(x_{1},\ldots,x_{n})=\prod_{i<j}(x_{i}-x_{j})\\
  J^{2}(x)=\prod_{i<j}(x_{i}-x_{j})^{2}=
  \prod_{\lambda\in\Delta_{+}(\mathfrak{a})}\lambda^{2}(x)
\end{gather*}
tandis que :
\begin{equation*}
  D(x)=\prod_{i<j}(x_{i}-x_{j})^{4}=
 \prod_{\lambda\in\Delta_{+}(\mathfrak{a})}\lambda^{2\mu(\lambda)}(x) 
\end{equation*}
où $\mu(\lambda)=2$ pour toute racine $\lambda$.

\subsection{Conclusion}
L'égalité :
\begin{equation*}
D(x)=\Vert\nabla p_{1}(x)\wedge\cdots\wedge\nabla p_{l}(x)\Vert^{2}
\end{equation*}
est vraie dans les deux cas des paires symétriques de rang maximal,
fausse dans les deux autres cas.

L'explication réside dans la comparaison entre
$\mathbb{R}[\mathfrak{p}]^{K}$ et l'algèbre
$\mathbb{R}[\mathfrak{g}]^{G}$ des fonctions polynômes $G$-invariantes
sur $\mathfrak{g}$. Tout d'abord, il est évident que la restriction à
$\mathfrak{p}$ d'une fonction $G$-invariante est une fonction
$K$-invariante. On dispose donc d'un homomorphisme d'algèbres :
\begin{equation*}
  Res_{\mathfrak{p}}:\mathbb{R}[\mathfrak{g}]^{G}\rightarrow
  \mathbb{R}[\mathfrak{p}]^{K}
\end{equation*}
dont on sait qu'il est surjectif, mis à part cas \flqq\;exceptionnels
\frqq. Quant à l'injectivité, la comparaison des degrés de
transcendance des deux algèbres montre qu'elle a lieu si et seulement
si la paire $(\mathfrak{g},\mathfrak{k})$ est de rang maximal. Au
total, $Res_{\mathfrak{p}}$ est un isomorphisme d'algèbres lorsque
$(\mathfrak{g},\mathfrak{k})$ est de rang maximal, et un système de
générateurs (homogènes, algébriquement indépendants)
$\tilde{p}_{1},\tilde{p}_{2},\ldots,\tilde{p}_{l}$ de
$\mathbb{R}[\mathfrak{g}]^{G}$ fournit, par restriction à
$\mathfrak{p}$, un système de générateurs $p_{1},p_{2},\ldots,p_{l}$
de $\mathbb{R}[\mathfrak{p}]^{K}$ ; les $d_{j}-1$ introduits plus haut
sont ce que l'on appelle les exposants de l'algèbre de Lie
$\mathfrak{g}$, et on a, pour tout $x$ dans $\mathfrak{p}$ :
\begin{equation*}
D(x)=\Vert\nabla p_{1}(x)\wedge\cdots\wedge\nabla p_{l}(x)\Vert^{2}=\phi(x)
\end{equation*}

On a rencontré dans ce texte deux fonctions polynômes sur
$\mathfrak{p}$ qui méritent la qualification de \flqq\;discriminant
\frqq~: elles sont $K$-invariantes, homogènes, sommes de carrés et
l'ensemble singulier $X$ de $\mathfrak{p}$ est exactement l'ensemble
des points de $\mathfrak{p}$ où s'annule l'une ou l'autre ce ces
fonctions. En comparant les restrictions à $\mathfrak{a}$ de ces
fonctions :
\begin{equation*}
  \phi(x)=\prod_{\lambda\in\Delta_{+}(\mathfrak{a})}\lambda^{2}(x)\mbox{,
  } D(x)=\prod_{\lambda\in\Delta_{+}(\mathfrak{a})}\lambda^{2\mu(\lambda)}(x)
\end{equation*}
on voit que $D$ est un multiple de $\phi$, qui semble donc être le
\flqq\;bon \frqq~ discriminant.

\end{document}